\documentclass{elsart}
\usepackage{amssymb}
\usepackage{amsmath}
\usepackage[english]{babel}

\newtheorem{theorem}{Theorem}

\newtheorem{lemma}[theorem]{Lemma}

\begin{document}
\begin{frontmatter}
\title{Fundamental solutions of homogeneous elliptic differential
operators.}
\author{Brice Camus}
\thanks[label2]{Work supported by the
\textit{SFB/TR12} project, \textit{Symmetries and Universality in
Mesoscopic Systems}.}
\address{Ruhr-Universit\"at Bochum, Fakult\"at f\"ur Mathematik,\newline%
Universit\"atsstr. 150, D-44780 Bochum, Germany.\newline Email :
brice.camus@univ-reims.fr}

\begin{abstract}
We compute fundamental solutions of homogeneous elliptic
differential operators, with constant coefficients, on
$\mathbb{R}^n$ by mean of analytic continuation of distributions.
The result obtained is valid in any dimension, for any degree and
can be extended to pseudodifferential operators of the same type.
\end{abstract}

\end{frontmatter}
Let be $P:=P(D_x)$ a pseudodifferential operator, with constant
coefficients, obtained by mathematical quantization of the
function $p$, i.e :
\begin{equation}
Pf(x)=\frac{1}{(2\pi)^n}\int\limits_{\mathbb{R}^{n}}
e^{i\left\langle \xi,x\right\rangle} p(\xi) \hat{f}(\xi)d\xi.
\end{equation}
Here, and in the following, the notation :
\begin{equation*}
\hat{f}(\xi)=\int\limits_{\mathbb{R}^n} e^{-i \left\langle y,
\xi\right\rangle} f(y) dy,
\end{equation*}
designs the Fourier transform of $f$. We note
$\mathcal{S}(\mathbb{R}^n)$ the Schwartz space and
$\mathcal{S}'(\mathbb{R}^n)$ the distributions on
$\mathcal{S}(\mathbb{R}^n)$. A fundamental solution for $P$ is a
distribution $\mathfrak{S}\in \mathcal{S}'(\mathbb{R}^n)$ such
that $P\mathfrak{S} =\delta$. Fundamental solutions play a major
role in the theory of PDE. For a large overview on this subject,
and applications, we refer to \cite{HOR1,HOR2}. Apart in some
trivial cases, there is few explicit characterizations of
fundamental solutions and most results concern the existence. The
case of order 3 homogeneous operators, in dimension 3, was treated
in \cite{Wag}. Always for $n=3$, the case of certain elliptic
homogeneous operators of degree 4 was also solved in \cite{Wag1}.
We are here interested in the case of a definite homogenous
polynomial $p_k$ on $\mathbb{R}^n$, i.e. $p_k(\xi)=0
\Leftrightarrow \xi=0$ and :
\begin{equation*}
p_k(\lambda \xi)=|\lambda|^{k} p_k(\xi).
\end{equation*}
Note that $k$ has to be even but we do not impose the spherical
symmetry of $p_k$. Strictly speaking, it is not necessary to
assume that $k\in\mathbb{N}$ and we can consider operators with a
conical singularity at the origin. The main motivation is that
such an operator can generalize the Laplacian (positivity and
ellipticity) but these operators also play a role in physical
optics. To state the main result, we introduce the spherical
average of $g\in\mathcal{S}(\mathbb{R}^n)$ w.r.t. the symbol $p$ :
\begin{equation}\label{spherical mean}
A(g)(r)=r^{n-1}\int\limits_{\mathbb{S}^{n-1}}
g(r\theta)p_k(\theta)^{-1}d\theta.
\end{equation}
Here $p_k(\theta)$ designs the restriction of $p_k$ on the sphere,
which defines a strictly positive function. Clearly, we have
$A(g)\in \mathcal{S}(\mathbb{R}_+)$ and we can extend $A(g)$ by 0
for $r<0$ to obtain a $L^2$ function on the line. The main result
is :
\begin{theorem}\label{main}
If $k<n$ a fundamental solution
$\mathfrak{S}\in\mathcal{S}'(\mathbb{R}^n)$ for $P$ is given by :
\begin{equation*}
\left\langle \mathfrak{S}, f
\right\rangle=\frac{1}{(2\pi)^n}\int\limits_{\mathbb{R}_{+}\times
\mathbb{S}^{n-1}} p_k(\theta)^{-1} \hat{f}(r\theta)
r^{n-1-k}drd\theta.
\end{equation*}
But, when $k\geq n$, we have :
\begin{equation*}
\left\langle \mathfrak{S}, f \right\rangle=-C_{k,n}
\frac{\partial^{k-1}A(\hat{f})}{\partial r^{k-1}}(0)+D_{k,n}
\int\limits_{\mathbb{R}_{+}}\log(u)\frac{\partial^{k}A(\hat{f})}{\partial
r^{k}}(u)du,
\end{equation*}
where $C_{k,n}$ and $D_{k,n}$ are universal constants given by
Eqs.(\ref{formula C_k},\ref{formula D_k}).
\end{theorem}
The reader can observe the analogy with the Laplacian, see e.g.
\cite{HOR1}. In particular, one has to distinguish the case of an
integrable (resp. non-integrable) singularity for $p_k(\xi)^{-1}$
for the $n$-dimensional Lebesgue measure.\medskip\\
\textit{Proof of Theorem \ref{main}.} Since $p_k\geq 0$, we define
a family $\mathfrak{p}(z)$ of distributions :
\begin{equation*}
\forall f\in\mathcal{S}(\mathbb{R}^n) : \left\langle
\mathfrak{p}(z),f\right\rangle= \frac{1}{(2\pi)^n}
\int\limits_{\mathbb{R}^n} p_k(\xi)^{z} \hat{f}(\xi)d\xi.
\end{equation*}
The r.h.s. is holomorphic for $\Re(z)>-1/k$ and, by continuity, we
obtain :
\begin{equation}\label{equ. of continuity}
\lim\limits_{z\rightarrow 0} \frac{1}{(2\pi)^n}
\int\limits_{\mathbb{R}^n} p_k(\xi)^{z}
\hat{f}(\xi)d\xi=f(0)=\left\langle \delta, f\right\rangle.
\end{equation}
The Laurent development of $\mathfrak{p}$ in $z=-1$ can be written
:
\begin{equation}\label{Laurent series}
\mathfrak{p}(z-1)=\sum\limits_{j=-1}^{-d} z^{j} \mu_j +\mu_0
+\sum\limits_{j=1}^\infty \mu_j z^j.
\end{equation}
This point is justified in Lemma \ref{analytic} below. But,
according to Eq.(\ref{equ. of continuity}), we have :
\begin{equation*}
\lim\limits_{z\rightarrow 0} \left\langle \mathfrak{p}(z-1),
P(D)f\right\rangle =\left\langle \delta, f\right\rangle.
\end{equation*}
It is then easy to check that $\mu_0$ is a fundamental solution
for $P$. Also, note that Eqs.(\ref{equ. of
continuity},\ref{Laurent series}) provide the set of non-trivial
relations :
\begin{equation}
P(D)\mu_j =0,\text{ } \forall j <0,
\end{equation}
in the sense of distributions of $\mathcal{S}'(\mathbb{R}^n)$. The
existence of such non-zero $\mu_j$, for $j<0$, implies the
non-uniqueness of fundamental solutions in
$\mathcal{S}'(\mathbb{R}^n)$.
\begin{lemma}\label{analytic}
The distributions $\mathfrak{p}(z-1)$ are meromorphic on
$\mathbb{C}$ with poles located at rational points
$z_{j,k}=-\frac{j}{k}$, $j\in\mathbb{N}$.
\end{lemma}
\textit{Proof.} Let $g=\hat{f}$. Using standard polar coordinates
we obtain :
\begin{equation*}
\left\langle \mathfrak{p}(z-1),
f\right\rangle=\frac{1}{(2\pi)^n}\int\limits_{\mathbb{R}_{+}\times\mathbb{S}^{n-1}}
(r^{k}p_k(\theta))^{z-1}g(r\theta) r^{n-1}drd\theta.
\end{equation*}
Here $p_k(\theta)$ is the restriction of $p_k$ to
$\mathbb{S}^{n-1}$. Next, if we define :
\begin{equation*}
y=(y_1,...,y_n)=(r p_k(\theta)^{\frac{1}{k}},\theta),
\end{equation*}
we obtain a very elementary formulation :
\begin{equation*}
\left\langle \mathfrak{p}(z-1),
f\right\rangle=\frac{1}{(2\pi)^n}\int\limits_{\mathbb{R}_{+}}
y_1^{k(z-1)} G(y_1)dy_1.
\end{equation*}
This new amplitude $G$ is obtained by pullback and integration :
\begin{equation}
G(y_1)= \int y^{*} (g(r,\theta)r^{n-1} |Jy|)dy_2 ...dy_n.
\end{equation}
We have $G\in \mathcal{S}(\mathbb{R}^{+})$ and
$G(y_1)=\mathcal{O}(y_1^{n-1})$ in $y_1=0$. Starting from :
\begin{equation*}
\frac{\partial^k}{\partial y_1^k}
y_1^{kz}=\prod\limits_{j=0}^{k-1} (kz-j) y_1^{k(z-1)},
\end{equation*}
after integrations by parts, we accordingly obtain that :
\begin{equation} \label{extended}
\left\langle \mathfrak{p}(z-1),
f\right\rangle=\frac{1}{(2\pi)^n}\prod\limits_{j=0}^{k-1}
\frac{1}{(kz-j)} \int\limits_{\mathbb{R}_{+}} y_1^{kz}
\partial^k_{y_1}G(y_1)dy_1.
\end{equation}
The integral of the r.h.s. is holomorphic in the strip
$\Re(z)>-\frac{1}{k}$. Finally, we can iterate the previous
construction to any order to obtain the result.
$\hfill{\blacksquare}$\medskip\\
Since $z=0$ is a simple pole, the constant term of the Laurent
series is :
\begin{equation*}
\frac{1}{(2\pi)^n k}\lim\limits_{z\rightarrow 0}
\partial_z( \prod\limits_{j=1}^{k-1} \frac{1}{(kz-j)}
\int\limits_{\mathbb{R}_{+}} y_1^{kz}
\partial^k_{y_1}G(y_1)dy_1).
\end{equation*}
For the calculations, we distinguish integrable and non-integrable singularities.\medskip\\
\textbf{Case of $k<n$}. The derivative of the rational function
provides :
\begin{equation}\label{formula C_k}
C_{k,n}=\frac{1}{(2\pi)^n k}(\partial_z( \prod\limits_{j=1}^{k-1}
\frac{1}{(kz-j)}))_{z=0} =\frac{(-1)^{k+1}}{(2\pi)^n} \left(
\frac{\gamma+\psi(k) }{\Gamma(k)} \right),
\end{equation}
where $\gamma$ is Euler's constant and $\psi$ the poly-gamma
function of order 0, i.e. :
\begin{gather*}
\gamma=\lim\limits_{m\rightarrow \infty} \sum\limits_{j=1}^{m} \frac{1}{k} -\log(m),\\
\psi(z)=\partial_z (\log(\Gamma(z))=\frac{\Gamma'(z)}{\Gamma(z)}.
\end{gather*}
Accordingly, since $G$ vanishes up to the order $n-1$ at the
origin, we have :
\begin{equation*}
C_{k,n} \int\limits_{\mathbb{R}_{+}} \partial^k_{y_1}G(y_1)dy_1=%
-C_{k,n} \partial^{k-1}_{y_1}G(0)=0.
\end{equation*}
On the other side, by derivation of the integral, we find the term
:
\begin{equation*}
\frac{(-1)^{k-1}}{(k-1)!} \int\limits_{\mathbb{R}_{+}}
\log(y_1)\partial^k_{y_1}G(y_1)dy_1.
\end{equation*}
Since $k<n$, we can integrate by parts to finally obtain :
\begin{equation*}
\int\limits_{\mathbb{R}_{+}} G(y_1)\frac{dy_1}{y_1^{k}}=
\int\limits_{\mathbb{R}_{+}\times \mathbb{S}^{n-1}}
p_k(\theta)^{-1} g(r\theta) r^{n-1-k}drd\theta.
\end{equation*}
Here, we have replaced the expression for our amplitude, via inversion of our diffeomorphism. This proves
the first statement of Theorem \ref{main}.$\hfill{\blacksquare}$\medskip\\
\textbf{Case of $k\geq n$.} The term attached to the derivative of
the rational function is non-zero and provides :
\begin{equation*}
-C_{k,n}
\partial^{k-1}_{y_1}G(0)+%
\frac{1}{(2\pi)^n} \frac{(-1)^{k-1}}{(k-1)!}
\int\limits_{\mathbb{R}_{+}}\log(y_1)\partial^k_{y_1}G(y_1)dy_1.
\end{equation*}
We compute first the derivative of $G$. Contrary to the case $k<n$
we cannot take the limit directly but we will reach the result
with the Schwartz kernel technic. By Fourier inversion formula we
have :
\begin{equation*}
\partial^{k-1}_{y_1}G(0)=\frac{1}{2\pi} \int\limits e^{-i\xi u}
(i\xi)^{k-1} G(u)dud\xi.
\end{equation*}
We can extend the integral w.r.t. $u$ on the whole line by
extending $G$ by zero for $u<0$. Going back to initial coordinates
provides :
\begin{equation*}
\partial^{k-1}_{y_1}G(0)
=\frac{1}{2\pi} \int\limits e^{-i\xi r} (i\xi)^{k-1}
r^{n-1}\int\limits_{\mathbb{S}^{n-1}} g(r\theta)
p_k(\theta)^{-1}d\theta dr d\xi.
\end{equation*}
The r.h.s. is exactly the derivative of order $k-1$ of $A(g)$
defined in Eq.(\ref{spherical mean}). By the same technic, the
logarithmic integral gives the non-local contribution :
\begin{equation*}
\frac{1}{(k-1)!} \int\limits_{\mathbb{R}_{+}} \log(u)
\frac{\partial^{k}A(g)}{\partial r^{k}}(u)du.
\end{equation*}
In particular the second universal constant of Theorem \ref{main}
is :
\begin{equation}\label{formula D_k}
D_{k,n}=\frac{1}{(2\pi)^n} \frac{(-1)^{k-1}}{(k-1)!}.
\end{equation}
This proves the second assertion of Theorem \ref{main}.
$\hfill{\blacksquare}$\medskip\\
The main result holds also holds for elliptic pseudo-differential
operators with homogeneous symbol of degree $\alpha\in
]0,\infty[$. Using $k=E(\alpha)+1$ integrations by parts in
Eq.(\ref{extended}), where $E(x)$ is the integer part of $x$, the
construction remains the same and all formulas can be analytically
continued w.r.t. the degree $\alpha>0$. This continuation is
trivial since the functions of Eqs.(\ref{formula C_k},\ref{formula
D_k}), defining respectively $C_{n,k}$ and $D_{n,k}$, are analytic
w.r.t. the variable $k>0$.

\end{document}